\documentclass{amsart}
\usepackage{amssymb}
\usepackage{amsmath}
\usepackage{textcomp}
\makeindex

\newcommand{\ba}{\begin{array}}
\newcommand{\eea}{\end{eqnarray}}
\newcommand{\ea}{\end{array}}

\newtheorem{definition}{Definition}[section]
\newtheorem{theorem}[definition]{Theorem}

\newtheorem{remark}[definition]{Remark}

\usepackage[matrix,arrow,curve]{xy}
\usepackage[utf8x]{inputenc}
\usepackage{tikz}

\begin{document}
\title[Rigidity of ADC Contact Structures]{Rigidity of ADC Contact Structures}
\author[S. Mukherjee]{Sauvik Mukherjee}

\maketitle
\section{Introduction} This paper is inspired by a question in \cite{Lazarev}. In \cite{Lazarev} Lazarev has proved the following result

\begin{theorem}(\cite{Lazarev})
\label{Lazarev}
All flexible weinstein fillings of a contact manifold $(Y,\xi)$ with $c_1(Y,\xi)=0$ have isomorphich integral cohomology.
\end{theorem}
{\bf Outline of the proof:} There exists a long exact sequence \[...\to H^{n-k}(W,\mathbb{Z})\to SH_k(W,\lambda)\to SH^+_k(W,\lambda)\to H^{n-k+1}(W,\mathbb{Z})\to...\] For flexible weinstein $W$, $SH_k(W,\lambda)=0$. So \[SH^+_k(W,\lambda)\cong H^{n-k+1}(W,\mathbb{Z})\] Now he introduces the concept of ADC (asymptotically dynamically convex) contact structure. He proved that if $W$ is flexible weinstein then $(Y,\xi)$ must be ADC moreover he proved that if $\partial W=Y$ is ADC then $SH^+_k(W,\lambda)$ depends only on $Y$ and not on $W$ and hence the proof.\\

Now he has asked the following question 

{\bf Question:} Is it true that any Liouville filling of $(Y,\xi)$ has isomorphich integral cohomologies if $(Y,\xi)$ is ADC?

In this paper we construct a counter example showing that the result is false. 

\begin{remark}
As we stated above that boundary of a flexible Weinstein domain must be ADC. Our counter example also shows that converse of this result is false.
\end{remark}

We end this section by stating some facts about ADC contact structures. Let $\gamma$ be a reeb orbit of $(Y,\xi)$ the degree of $\gamma$ is \[|\gamma|=\mu_{cz}(\gamma)+n-3\] where $\mu_{cz}(\gamma)$ is the Conley-Zehnder index of $\gamma$. Now we define \[A(\gamma)=\int_{S^1}\gamma^*\alpha\] where $\alpha$ is a contact form for $\xi$. Set $\mathcal{P}^{<D}(Y,\alpha)$ to be the set of contractible reeb orbit of $(Y,\alpha)$ such that $A(\gamma)<D$. We say $\alpha_1\geq \alpha_2$ for two contact forms $\alpha_i,\ i=1,2$ if $\alpha_1=f\alpha_2$ and $f\geq1$. 

\begin{definition}(\cite{Lazarev})
\label{ADC}
$(Y,\xi)$ is said to be ADC(asymptotically dynamically convex) if there exist contact forms $\alpha_1\geq \alpha_2\geq \alpha_3\geq....$ for $\xi$ and positive numbers $D_1<D_2<D_3<...$ going to infinity such that all elements of $\mathcal{P}^{<D_k}(Y,\alpha_k)$ have positive degree.
\end{definition}
Define the following 
\[\mathcal{W}^{<D}(\Lambda,Y,\alpha)=\{[w],\ w=c_1c_2..c_k:A(w)=\Sigma A(c_i)<D\}\] where $c$ is a reeb chord of the Legendrian $\Lambda$ and \[A(c)=\int_c\alpha\]  $w$ is zero in $\pi_1(Y,\alpha)$ where $w$ as an element of  $\pi_1(Y,\alpha)$ by inserting paths in $\lambda$ connecting one end of $c_1$ with the base point and connecting the end point of $c_i$ with the end point of $c_{i+1}$ for all $i$. The equivalence class $[w]$ is defined as $c_1c_2..c_n$ is equivalent to $c_2c_3...c_nc_1$.

\begin{theorem}(\cite{Lazarev})
\label{RC}
Let $\Lambda^{n-1}\subset (Y^{2n-1}_-,\alpha_-),\ n\geq 3$ be Legendrian sphere. For any $D>0$ there exists $\epsilon=\epsilon(D)>0$ such that if $(Y_+,\alpha_+)$ is the result of contact surgery on $U^{\epsilon}(\Lambda,\alpha_-)$ then there is grading preserving bijection between $\mathcal{P}^{<D}(Y_+,\alpha_+)$ and $\mathcal{P}^{<D}(Y_-\alpha_-)\cup \mathcal{W}^{<D}(\Lambda,Y_-,\alpha_-)$: if $\gamma_w$ is the orbit corresponding to the word of chords $w=c_1...c_k$ then \[|\gamma_w|=|w|+n-3=\Sigma |c_i|+n-3\] where $|c|=\mu_{cz}(c)-1$.
\end{theorem}

\section{Weinstein Handles} \label{WH} In this section we present an alternative elementary method of Weinstein handle attachments for the case of critical index in dimension $2n=8$. It is different from \cite{Weinstein} as it allows us to give explicit formula for the attaching map but is not applicable to all Legendrian spheres.\\

Consider the following spaces \[M=\{(p,q)\in \mathbb{R}^4\times \mathbb{R}^4:|p|^2=\epsilon_M\}\] \[N=\{(p,q)\in \mathbb{R}^4\times \mathbb{R}^4:|q|^2=\epsilon_N\}\]

Consider coordinates on $M\cap \{p_4>0\}$ as $(p_1,p_2,p_3,p_4,q_1,q_2,q_3,q_4)\to(p_1,p_2,p_3,q_1,q_2,q_3,q_4)$ and on $N\cap \{q_4>0\}$ as $(p_1,p_2,p_3,p_4,q_1,q_2,q_3,q_4)\to(p_1,p_2,p_3,p_4,q_1,q_2,q_3)$. On $M\cap\{p_4>0\}$ and on $N\cap\{q_4>0\}$ set respectively \[p'_4=(\epsilon_M-\Sigma_1^3p_i^2)^{1/2},\ q_4=(\epsilon_N-\Sigma_1^3q_i^2)\]

Consider the one forms as follows \[(\alpha_4)_{|M}=2\Sigma_1^3(q_i-q_4\frac{p_i}{p'_4})dp_i+\Sigma_1^3p_idq_i+p'_4dq_4\] \[(\alpha_4)_{|N}=\Sigma_1^3(p_i-p_4\frac{q_i}{q'_4})dq_i+2\Sigma_1^3q_idp_i+2q'_4dp_4\]

Now we wish to construct a strict contactomorphism from $M$ to $N$ i.e, \[F:M\to N,\ F(p_1,p_2,p_3,q_1,q_2,q_3,q_4)=(r_1,r_2,r_3,r_4,\frac{\epsilon_N}{\epsilon_M}p_1,\frac{\epsilon_N}{\epsilon_M}p_2,\frac{\epsilon_N}{\epsilon_M}p_3)\] such that $F^*(\alpha_4)_N=(\alpha_4)_M$ with the condition 
\begin{equation}
\label{LC}
F(p_1,p_2,p_3,0,0,0,0)=(r'_1,r'_2,r'_3,r'_4,\frac{\epsilon_N}{\epsilon_M}p_1,\frac{\epsilon_N}{\epsilon_M}p_2,\frac{\epsilon_N}{\epsilon_M}p_3)=Given\ Legendrian
\end{equation}

So we write \[(\alpha_4)_N=2\frac{\epsilon_N}{\epsilon_M}\Sigma_1^3p_idr_i+2\frac{\epsilon_N}{\epsilon_M}p'_4dr_4+\frac{\epsilon_N}{\epsilon_M}\Sigma_1^3(r_i-r_4\frac{p_i}{p'_4})dp_i=(\alpha_4)_M\] By substituting $dr_j=\Sigma_1^3 \frac{\partial r_j}{\partial p_i}dp_i+\Sigma_1^4\frac{\partial r_j}{\partial q_i}dq_i$ in the above we get 

From now we write $p'_4(p_1,p_2,p_3)=p_4(p_1,p_2,p_3)$ for convenient notation.
 
\begin{equation}
\label{K1}
2\frac{\epsilon_N}{\epsilon_M}\Sigma p_i\frac{\partial r_i}{\partial q_j}=p_j,\ j=1,2,3,4
\end{equation}
Integrating we get 
\begin{equation}
\label{K11}
\frac{\epsilon_N}{\epsilon_M}\Sigma p_ir_i=\frac{1}{2}p_jq_j+A_j,\ j=1,2,3,4
\end{equation}
where $A_j$ is a function which does not depend on $q_j$ for $j=1,2,3,4$. Similarly we get 
\begin{equation}
\label{K2}
2\frac{\epsilon_N}{\epsilon_M}\Sigma p_i\frac{\partial r_i}{\partial p_j}+\frac{\epsilon_N}{\epsilon_M}(r_j-r_4\frac{p_j}{p_4})=2(q_j-q_4\frac{p_j}{p_4}),\ j=1,2,3
\end{equation}

Differentiating equation \ref{K11} with respect to $p_j$ for $j=1,2,3$ we get 

\begin{equation}
\label{K3}
\frac{\epsilon_N}{\epsilon_M}\Sigma p_i\frac{\partial r_i}{\partial p_j}+\frac{\epsilon_N}{\epsilon_M}(r_j-r_4\frac{p_j}{p_4})=\frac{q_j}{2}+\frac{\partial A_j}{\partial p_j}
\end{equation}

Comparing equation \ref{K2} and \ref{K3} we get 
\begin{equation}
\label{K4}
\frac{\epsilon_N}{\epsilon_M}(r_j-r_4\frac{p_j}{p_4})=q_j-2(q_j-q_4\frac{p_j}{p_4})+2\frac{\partial A_j}{\partial p_j}
\end{equation}

Now from condition \ref{LC} and equation \ref{K11} we get $\frac{\epsilon_N}{\epsilon_M}\Sigma p_ir'_i=A_j(p_1,p_2,p_3,0,0,0,0)=A'_j(say)$. So $A_j=A'_j+B_j$ such that $B_j$ does not depend on $q_j$ and \[B_j(p_1,p_2,p_3,0,0,0)=0\] So $A'_1=A'_2=A'_3=A'_4$. So we have \[\frac{\epsilon_N}{\epsilon_M}\Sigma p_ir_i=\frac{1}{2}p_jq_j+A'_j+B_j\] which gives us the system of equations 

\begin{equation}
\label{K5}
B_k-B_l=\frac{1}{2}(p_lq_l-p_kq_k)+A'_l-A'_k=\frac{1}{2}(p_lq_l-p_kq_k)
\end{equation}
whose coefficient matrix has determinant zero. The solutions are \[B_j=B_4+\frac{1}{2}(p_4q_4-p_jq_j),\ j\neq 4\] we choose $B_4$ such that $B_j$ does not depend on $q_j$. So take \[B_4=\frac{1}{2}\Sigma_{i\neq 4}p_iq_i\] Now we just need to solve system of equations consisting of equation \ref{K11} for $j=4$ and equations \ref{K4} with the above substitution of $B_j$ hence of $A_j$. This gives \[\frac{\epsilon_N}{\epsilon_M}r_4=\frac{p_4}{\epsilon_M}[\frac{1}{2}\Sigma p_iq_i+\Sigma_{i\neq 4}p_iq_i-2\frac{q_4}{p_4}\Sigma_{i\neq 4}p_i^2+A'_4-2\Sigma_{i\neq 4}p_i\frac{\partial A'_i}{\partial p_j}]\] Now the rest of $r_i$'s are found using equation \ref{K4}.

\begin{remark}
\label{loose}
The given Legendrian in the above construction can not be loose. To see this first observe that $N\cap \{q_4>0\}$ and $N\cap \{q_4<0\}$ covers $N$ except the $\{q_4=0\}$. Without loss of generality consider the Darboux chart on $N\cap \{q_4>0\}$ as \[\phi(p_1,p_2,p_3,p_4,q_1,q_2,q_3)=(p_4+\Sigma_1^3\frac{q_i}{{q'}_4}p_i,e_1,e_2,e_3,q_1,q_2,q_3)\] where $e_j=-(\frac{p_j}{2{q'}_4}+\frac{q_j}{{q'}_4^3}\Sigma_1^3p_iq_i+p_4\frac{q_j}{2{q'}_4^2})$. From \ref{LC} we write \[\phi(F(p_1,p_2,p_3,0,0,0,0))=\phi(r'_1,r'_2,r'_3,r'_4,\frac{\epsilon_N}{\epsilon_M}p_1,\frac{\epsilon_N}{\epsilon_M}p_2,\frac{\epsilon_N}{\epsilon_M}p_3)\] which is equal to $(r'_4+\Sigma_1^3\frac{p_i}{p'_4}r'_i,E_1,E_2,E_3,\frac{\epsilon_N}{\epsilon_M}p_1,\frac{\epsilon_N}{\epsilon_M}p_2,\frac{\epsilon_N}{\epsilon_M}p_3)$ where $E_j=e_j(r'_1,r'_2,r'_3,r'_4,\frac{\epsilon_N}{\epsilon_M}p_1,\frac{\epsilon_N}{\epsilon_M}p_2,\frac{\epsilon_N}{\epsilon_M}p_3)$ So $\phi(F(p_1,p_2,p_3,0,0,0,0))$ can not have a cusp in the front projection. As loose legendrians has double cusp one of the cusp must be in $N\cap \{q_4>0\}$ which gives a contradictiion.
\end{remark}

\section{The Counter Example}
In this section we construct the required counter example, i.e, Weinstein domains $W_1$ and $W_2$ such that the boundaries of $W_i$'s are same and the contact structure on the boundary is ADC(asymptotically dynamically convex). First we set $W_1$ to be the space by attaching standard Weinstein handle $H_W$ to $D^8$ (eight dimensional disc) along some unknot in $\partial D^8=S^7$. Note that $H_W$ has boundary as \[\{(p,q)\in \mathbb{R}^4\times \mathbb{R}^4:|p|^2=2: |q|\leq \epsilon_N\}\cup \{(p,q)\in \mathbb{R}^4\times \mathbb{R}^4:|q|^2=\epsilon_N: |p|\leq 2\}\] Let $N'=\{(p,q)\in \mathbb{R}^4\times \mathbb{R}^4:|q|^2=\epsilon_N: |p|\leq 2\}\subset \partial H_W$. Now we construct $W_2$. For this first consider $F:M\to N'$ as in the section \ref{WH} with $r'_i=0$. Observe that $\{(0,0,0,0,q_1,q_2,q_3)\subset N'\}$ is a Legendrian in $N'$. \\

We consider a $8$-dimensional weinstein handle $H$ with boundary \[\partial H=M\cup \bar{N}\] where $\bar{N}$ will be explained shortly. So $F$ attaches $H$ with $H_W$ and hence with $W_1$. Now we consider another weinstein handle $\bar{H}\subset \mathbb{R}^4\times \mathbb{R}^4$ with boundary \[\partial \bar{H}=\bar{M}\cup F(M)\] where $\bar{M}$ will also be explained shortly.\\

We shall define $\bar{M}$ and $\bar{N}$ that will help us construct an attaching map between $\bar{M}$ and $\bar{N}$ so that we shall be able to attach $H$ and $\bar{H}$. This will be the required $W_2$. Observe that boundary of $W_i$'s are same.\\

First we define $\bar{N}$. Set $\epsilon_M=\epsilon_N$ as we have chosen the handle $H$. Then \[\bar{N}=\{(P,Q)\in \mathbb{R}^4\times \mathbb{R}^4\}\] such that near the corner $P_i(p,q)=r_i(p,q),\ Q_i(p,q)=p_i$ and $P$ and $Q$ are such that $\bar{N}$ is transversal to the Liouville vectorfield. Now we define $\bar{M}$ similarly \[\bar{M}=\{(P,Q)\in \mathbb{R}^4\times \mathbb{R}^4\}\] such that near the corner $P_i(p,q)=r_i(p,q),\ Q_i(p,q)=p_i$ and $P$ and $Q$ are such that $\bar{N}$ is transversal to the Liouville vectorfield.Here the pull back condition does not involve derivatives as we only need to match the coefficients of $dp_i$ and $dq_i$ which are $2Q_i$ and $P_i$ respectively.\\

Now we show that $\partial W_i$ is ADC. We just need to show $\partial W_1$ is ADC as $\partial W_2=\partial W_1$. As we have attached $H_W$ along an unknot we just need to consider $a^k$ in $\mathcal{W}^{<D}(\Lambda,Y_-,\alpha_-)$ as in \ref{RC} where $Y_-$ is $S^7$, $\Lambda$ is the unknot and $a$ is the reeb chord in the unknot. If $\gamma_w$ be the reeb orbit corresponding to $a^k$ then $|\gamma_w|=3k+1>0$.

\begin{remark}
By remark \ref{loose} above a construction similar to ours can not produce counter example to \ref{Lazarev}.
\end{remark}

\end{document}